\pdfoutput=1
\RequirePackage{ifpdf}
\ifpdf 
\documentclass[pdftex]{sigma}
\else
\documentclass{sigma}
\fi

\numberwithin{equation}{section}
\newtheorem{Theorem}{Theorem}[section]
\newtheorem{Corollary}[Theorem]{Corollary}
{ \theoremstyle{definition}
\newtheorem{Remark}[Theorem]{Remark} }

\begin{document}

\newcommand{\arXivNumber}{1402.1569}

\allowdisplaybreaks

\renewcommand{\PaperNumber}{103}

\FirstPageHeading

\ShortArticleName{On Certain Wronskians of Multiple Orthogonal Polynomials}

\ArticleName{On Certain Wronskians of Multiple Orthogonal\\
Polynomials}

\Author{Lun ZHANG~$^\dag$ and Galina FILIPUK~$^\ddag$}

\AuthorNameForHeading{L.~Zhang and G.~Filipuk}

\Address{$^\dag$~School of Mathematical Sciences and Shanghai Key Laboratory for Contemporary\\
\hphantom{$^\dag$}~Applied Mathematics, Fudan University, Shanghai 200433, People's Republic of China}
\EmailD{\href{mailto:lunzhang@fudan.edu.cn}{lunzhang@fudan.edu.cn}}
\URLaddressD{\url{http://homepage.fudan.edu.cn/lunzhang/}}

\Address{$^\ddag$~Faculty of Mathematics, Informatics and Mechanics, University of Warsaw,\\
\hphantom{$^\ddag$}~Banacha 2, Warsaw, 02-097, Poland}
\EmailD{\href{mailto:filipuk@mimuw.edu.pl}{filipuk@mimuw.edu.pl}}
\URLaddressD{\url{http://www.mimuw.edu.pl/~filipuk/}}

\ArticleDates{Received August 01, 2014, in f\/inal form October 27, 2014; Published online November 04, 2014}

\Abstract{We consider determinants of Wronskian type whose entries are multiple ortho\-go\-nal polynomials associated with
a~path connecting two multi-indices.
By assuming that the weight functions form an algebraic Chebyshev (AT) system, we show that the polynomials represented
by the Wronskians keep a~constant sign in some cases, while in some other cases oscillatory behavior appears, which
generalizes classical results for orthogonal polynomials due to Karlin and Szeg\H{o}.
There are two applications of our results.
The f\/irst application arises from the observation that the~$m$-th moment of the average characteristic polyno\-mials for
multiple orthogonal polynomial ensembles can be expressed as a~Wronskian of the type~II multiple orthogonal polynomials.
Hence, it is straightforward to obtain the distinct be\-havior of the moments for odd and even~$m$ in a~special multiple
orthogonal ensemble~-- the AT ensemble.
As the second application, we derive some Tur\'an type inequalities for multiple Hermite and multiple Laguerre
polynomials (of two kinds).
Finally, we study numerically the geometric conf\/iguration of zeros for the Wronskians of these multiple orthogonal
polynomials.
We observe that the zeros have regular conf\/igurations in the complex plane, which might be of independent interest.}


\Classification{05E35; 11C20; 12D10; 26D05; 41A50}

\section{Introduction and statement of the main results}

\subsection{Determinants whose entries are orthogonal polynomials}
In a~classical paper~\cite{KS}, Karlin and Szeg\H{o} developed an interesting and general theory regarding the
determinants whose entries are orthogonal polynomials~\cite{Chi,Gau,Sze}.
They showed that the polynomials represented by certain determinants whose elements are orthogonal polynomials keep
a~constant sign in some cases, while in some other cases the polynomials are oscillatory.
More precisely, let
\begin{gather*}
Q_n(x)=k_n(-x)^n+\cdots,
\qquad
k_n>0,
\qquad
n\in\mathbb{N} = \{0,1,2,3,\ldots \},
\end{gather*}
be a~sequence of orthogonal polynomials with respect to an arbitrary measure whose distribution function has an inf\/inite
number of increasing points.
The Wronskian of these polynomials is then def\/ined~by
\begin{gather*}
W(n,l;x) :=W\left(Q_n(x), Q_{n+1}(x), \dots, Q_{n+l-1}(x)\right)
\\
\phantom{W(n,l;x)}
 = \det \left(
\begin{matrix}
Q_{n}(x) & Q_{n+1}(x) & \cdots & Q_{n+l-1}(x)
\\
Q_{n}'(x) & Q_{n+1}'(x) & \cdots & Q_{n+l-1}'(x)
\\
\vdots & \vdots & \vdots & \vdots
\\
Q_{n}^{(l-1)}(x) & Q_{n+1}^{(l-1)}(x) & \cdots & Q_{n+l-1}^{(l-1)}(x)
\\
\end{matrix}
\right).
\end{gather*}
By~\cite[Theorems 1 and 2]{KS}, it is known that, for~$l$ even,
\begin{gather*}
(-1)^{l/2}W(n,l;x)>0,
\qquad
x\in\mathbb{R},
\end{gather*}
i.e., the Wronskian keeps a~constant sign for all real~$x$; if~$l$ is odd, then $W(n,l;x)$ has exactly~$n$ simple real
zeros and the zeros of $W(n,l;x)$ and $W(n+1,l;x)$ strictly interlace.

Another important class of determinants considered in~\cite{KS} is the Hankel determinant
\begin{gather*}
T(n,l;x) :=T\left(Q_n(x), Q_{n+1}(x), \ldots, Q_{n+l-1}(x)\right)
\\
\phantom{T(n,l;x)}
 = \det \left(
\begin{matrix}
Q_{n}(x) & Q_{n+1}(x) & \cdots & Q_{n+l-1}(x)
\\
Q_{n+1}(x) & Q_{n+2}(x) & \cdots & Q_{n+l}(x)
\\
\vdots & \vdots & \vdots & \vdots
\\
Q_{n+l-1}(x) & Q_{n+l}(x) & \cdots & Q_{n+2l-2}(x)
\\
\end{matrix}
\right),
\end{gather*}
which is called the Tur\'{a}nian.
Karlin and Szeg\H{o} showed that, if~$l$ is even, $T(n,l;x)$ has the sign $(-1)^{l/2}$ on the interval~$I$ for the
following three classical systems of orthogonal polynomials~\cite{Sze}:
\begin{itemize}\itemsep=0pt
\item $Q_n(x)=P_n^{(\lambda)}(x)/P_n^{(\lambda)}(1)$, $\lambda > -1/2$ and $I=(-1,1)$, where $P_n^{(\lambda)}(x)$ are
the ultraspherical polynomials,
\item $Q_n(x)=L_n^{(\alpha)}(x)/L_n^{(\alpha)}(0)$, $\alpha > -1$ and $I=(0,+\infty)$, where $L_n^{(\alpha)}(x)$ are the
Laguerre polynomials,
\item $Q_n(x)=H_n(x)$ and $I=(-\infty,+\infty)$, where $H_n(x)$ are the Hermite polynomials.
\end{itemize}
The strategy of proofs is to represent the Hankel determinants in terms of the Wronskian of orthogonal polynomials of
another type.
Note that, if $l=2$, one has
\begin{gather*}
T(n,2;x)=Q_{n+1}^2(x)-Q_{n}(x)Q_{n+2}(x)>0.
\end{gather*}
This inequality is called the Tur\'{a}n inequality, which was f\/irst proved for the Legendre polynomials
$P_{n}(x)=P_n^{(1/2)}(x)$~\cite{Szego1,Turan} and inspired the work of Karlin and Szeg\H{o}.
The analogous results for determinants involving orthogonal polynomials associated with discrete weights are also
presented in~\cite{KS}.

Nowadays, determinants whose entries are orthogonal polynomials still attract much attention.
For instance, the Wronskians of orthogonal polynomials also appear in random matrix theory; cf.~\cite{BH2001,MN} and
Section~\ref{sec:moments} below.
The relationship between the Wronskian of orthogonal polynomials and the Hankel determinant of polynomials is clarif\/ied
by Leclerc in~\cite{Leclerc}, and further generalized by Dur\'{a}n~\cite{Duran}.
In addition, it comes out that Tur\'an inequality holds not only for a~large class of orthogonal polynomials including
the most classical orthogonal polynomials (cf.~\cite{Berg,Busto,Elbert2,Elbert1,Gas,Gas2,Krasikov,Szasz1,Szasz2}), but
also for many special functions and their~$q$-analogues with important applications; we refer
to~\cite{ABusto,Baricz07,BariczIsmail,BJP,Baricz,CNV,IsmailL,Laforgia,Lorch,Nuttall,Sko} and the re\-fe\-ren\-ces therein for
the development of that aspect.
Other studies of the Wronkians of orthogonal polynomials can be found in~\cite{Mourad,Samuel,Robert}.

In this paper, we are concerned with the Wronskians of multiple orthogonal polynomials.
Since multiple orthogonal polynomials are generalizations of orthogonal polynomials, our results extend the
aforementioned results for orthogonal polynomials.
In what follows, we f\/irst give a~brief introduction to multiple orthogonal polynomials and f\/ix the notations used
throughout this paper, and next state the main results and outline the rest of the paper.

\subsection{Multiple orthogonal polynomials and algebraic Chebyshev (AT) systems}

Multiple orthogonal polynomials are polynomials of one variable which are def\/ined by orthogonality relations with
respect to~$r$ dif\/ferent weights $w_1,w_2,  \ldots,w_r$,
where $r \geq 1$ is a~positive integer.
They originated from Hermite--Pad\'e approximation in the context of irrationality and transcendence proofs in number
theory, and they were further developed in approximation theo\-ry; cf.~\cite{Apt98,Apt99,Bust,Gon,NikSor} and
surveys~\cite{Apt,WVA99,WVA06}.

Let $\vec{n} = (n_1,n_2,\ldots,n_r) \in \mathbb{N}^r$ be a~multi-index of size $|\vec{n}| = n_1+n_2+\dots+n_r$ and
suppose $w_1,w_2,  \ldots,w_r$
are~$r$ weights with supports on the real axis.
There are two types of multiple orthogonal polynomials.
The \emph{type I multiple orthogonal polynomials} are given by the vector $(A_{\vec{n},1},\ldots, A_{\vec{n},r})$, where
$A_{\vec{n},j}$ is a~polynomial of degree $\leq n_j-1$, such that the linear form
\begin{gather}
\label{def:Qn}
Q_{\vec{n}}(x) = \sum\limits_{j=1}^r A_{\vec{n},j}(x)w_j(x)
\end{gather}
satisf\/ies
\begin{gather}
\label{eq:1.6}
\int Q_{\vec{n}}(x) x^k \, dx = 0,
\qquad
k=0,1,\ldots,|\vec{n}|-2.
\end{gather}
By setting the normalization condition
\begin{gather}
\label{eq:1.7}
\int Q_{\vec{n}}(x) x^{|\vec{n}|-1} \, dx = 1,
\end{gather}
the equations~\eqref{eq:1.6}, \eqref{eq:1.7} form a~linear system of $|\vec{n}|$ equations for the unknown coef\/f\/icients
of $A_{\vec{n},1},\ldots,A_{\vec{n},r}$.
The multi-index $\vec{n}$ is called {\it normal} if this linear system has a~unique solution, i.e., the polynomials of
vector $(A_{\vec{n},1},\ldots,A_{\vec{n},r})$ exist uniquely.
The \emph{type II multiple orthogonal polynomial} is the monic polynomial $P_{\vec{n}}(x) = x^{|\vec{n}|} + \cdots$ of
degree $|\vec{n}|$ satisfying the conditions
\begin{gather}
\int P_{\vec{n}}(x) x^k w_1(x)\, dx  = 0,
\qquad
k=0,1,\ldots,n_1-1,
\nonumber
\\
\cdots\cdots\cdots\cdots\cdots\cdots\cdots\cdots\cdots\cdots\cdots\cdots\cdots\cdots\cdots\cdots\label{eq:orthogonality of Pn}
\\
\int P_{\vec{n}}(x) x^k\, w_r(x)\, dx  = 0,
\qquad
k=0,1,\ldots,n_r-1.
\nonumber
\end{gather}
The polynomials $P_{\vec{n}}$ exist and are unique whenever $\vec{n}$ is a~normal index.

Under certain additional conditions on~$r$ weights, we can ensure the uniqueness and existence of multiple orthogonal
polynomials.
One of such conditions is that the weight functions form a~so-called algebraic Chebyshev (AT) system; cf.~\cite[Section 23.1.2]{Ismail}.

A Chebyshev system $\{\varphi_1,\ldots,\varphi_n\}$ on $[a,b]$ is a~system of~$n$ linearly independent functions such
that every linear combination $\sum\limits_{k=1}^n a_k \varphi_k$ has at most $n-1$ zeros on $[a,b]$.
Equivalently, this means that
\begin{gather*}
\det
\begin{pmatrix}
\varphi_1(x_1) & \varphi_1(x_2) & \cdots & \varphi_1(x_n)
\\
\varphi_2(x_1) & \varphi_2(x_2) & \cdots & \varphi_2(x_n)
\\
\vdots & \vdots & \cdots & \vdots
\\
\varphi_n(x_1) & \varphi_n(x_2) & \cdots & \varphi_n(x_n)
\end{pmatrix}
\neq 0,
\end{gather*}
for every choice of~$n$ dif\/ferent points $x_1,\ldots,x_n \in [a,b]$.
To see this, suppose $x_1, \ldots, x_n$ are such that the determinant is zero, then there is a~linear combination of the
rows $\sum\limits_{k=1}^n b_k \varphi_k$ that gives a~zero row.
We then obtain a~linear combination of functions $\varphi_k$ admitting~$n$ zeros at $x_1,\ldots,x_n$, which is
a~contradiction.

A system of~$r$ weights $(w_1,\ldots,w_r)$ is an AT system for the multi-index $\vec{n}$ if each $w_j$ is def\/ined on
a~f\/ixed interval $[a,b]\subseteq \mathbb{R}$ such that
\begin{gather*}
\big\{w_1,xw_1,\ldots,x^{n_1-1}w_1,w_2,xw_2,\ldots,x^{n_2-1}w_2,\ldots, w_r,xw_r,\ldots,x^{n_r-1}w_r\big\}
\end{gather*}
is a~Chebyshev system on $[a, b]$.
If $\vec{n}$ is a~multi-index such that the weights $(w_1,\ldots,w_r)$ form an AT system for every index $\vec{m}$
satisfying $\vec{m}\leq \vec{n}$ (in the componentwise sense, that is, $m_j\leq n_j$, $j=1,\ldots,r$),
by~\cite[Theorem~23.2]{Ismail}. We have that $\vec{n}$ is a~normal index, which implies the existence and uniqueness
of the polynomials $P_{\vec{n}}$.
In particular, the weights for many classical multiple orthogonal polynomials (including multiple Hermite polynomials,
multiple Laguerre polynomials, Jacobi--Pi\~{n}eiro polynomials) belong to the AT systems.

For more information about multiple orthogonal polynomials, we refer to Aptekarev et al.~\cite{Apt,AptBraWVA},
Coussement and Van Assche~\cite{WVAEC}, Nikishin
and Sorokin~\cite[Chapter 4, \S~3]{NikSor}, Ismail~\cite[Chapter 23]{Ismail}
and Filipuk, Van Assche and Zhang~\cite{MOP}.

\subsection{Statement of the main results}
\label{subsec:main results}
To state our main results, we need to def\/ine the Wronskian of multiple orthogonal polynomials.
Let us consider a~sequence of multi-indices $(\vec{n}_0,\vec{n}_1,\ldots,\vec{n}_{l-1})$, $l\in\mathbb{Z}^+=
\{1,2,3,\ldots \}$, such that
\begin{itemize}\itemsep=0pt
\item $\vec{n}_0=\vec{n}$ for a~given initial multi-index $\vec{n}$,
\item $ |\vec{n}_j|=|\vec{n}|+j$ for $j=1,\ldots, l-1$,
\item $ \vec{n}_0\leq \vec{n}_1\leq \dots \leq \vec{n}_{l-2}\leq\vec{n}_{l-1}$ (componentwise).
\end{itemize}
Therefore, $(\vec{n}_0,\vec{n}_1,\ldots,\vec{n}_{l-1})$ def\/ines a~path connecting $\vec{n}$ to $\vec{n}_{l-1}$, where in
each step the multi-index $\vec{n}_k$ is increased by one in exactly one direction.

For every such kind of a~f\/ixed path, we def\/ine the associated Wronskian of multiple orthogonal polynomials~by
\begin{gather}
W(\vec{n},l;x) :=W\left(P_{\vec{n}_0}(x), P_{\vec{n}_{1}}(x), \ldots, P_{\vec{n}_{l-1}}(x)\right)
\nonumber
\\
\phantom{W(\vec{n},l;x)}
 = \det \left(
\begin{matrix}
P_{\vec{n}_0}(x) & P_{\vec{n}_{1}}(x) & \cdots & P_{\vec{n}_{l-1}}(x)
\\
P_{\vec{n}_0}'(x) & P_{\vec{n}_{1}}'(x) & \cdots & P_{\vec{n}_{l-1}}'(x)
\\
\vdots & \vdots & \vdots & \vdots
\\
P_{\vec{n}_0}^{(l-1)}(x) & P_{\vec{n}_{1}}^{(l-1)}(x) & \cdots & P_{\vec{n}_{l-1}}^{(l-1)}(x)
\end{matrix}
\right),
\label{def:wronskian}
\end{gather}
where $P_{\vec{n}}$ is the type II multiple orthogonal polynomial given in~\eqref{eq:orthogonality of Pn}.
Clearly, $W(\vec{n},l;x)$ is a~polynomial in~$x$ depending on the parameters $\vec{n}$, $l$ and the path starting from
$\vec{n}$.
We shall also use the notation $W(\vec{a},\vec{b},\ldots,\vec{e};x)$ to emphasize the dependence on a~specif\/ic path
consisting of the multi-indices $(\vec{a},\vec{b},\ldots,\vec{e})$ if necessary.

Our main results are stated as follows.

\begin{Theorem}
\label{thm:l even}
Suppose that the weights $(w_1,w_2,\ldots,w_r)$ form an AT system on $[a,b]$ for all the multi-indices $\vec
{n}\in\mathbb{N}^r$, then we have
\begin{gather}
\label{eq: l even}
W(\vec{n},l;x)>0,
\qquad
x\in\mathbb{R},
\end{gather}
if~$l$ is even, where $W(\vec{n},l;x)$ is defined in~\eqref{def:wronskian}.
\end{Theorem}

Note that our assumption on the weights ensures the existence of multiple orthogonal polynomials (see the discussion at
the end of the previous section), thus the Wronskian is well-def\/ined.
If~$l$ is odd, then we have the following result.

\begin{Theorem}
\label{thm:l odd}
Let $w_1,w_2,\ldots,w_r$ be~$r$ weights as in Theorem~{\rm \ref{thm:l even}} and let~$l$ be odd.
For each fixed multi-index $\vec{n}$ the polynomials $ W(\vec{n},l;x) $ have exactly $|\vec{n}|$ simple zeros on the
real axis.
Furthermore, given two paths consisting of~$l$ multi-indices such that the last $l-1$ multi-indices of one path starting
from $\vec{n}$ coincide with the first $l-1$ multi-indices of the other path ending at~$\vec{m}$ $($which also means
$|\vec{m}|=|\vec{n}|+l$ and $\vec{n}\leq\vec{m})$, then the real zeros of two associated Wronskians strictly interlace.
\end{Theorem}

In case $l=1$, this theorem states that the type II multiple orthogonal polynomial $P_{\vec{n}}$ whose weights form an
AT system has $|\vec{n}|$ zeros and the zeros of $P_{\vec{n}}$ and $P_{\vec{n}+\vec{e}_k}$ ($k=1,\ldots,r$) interlace,
where $\vec{e}_k = (0,\ldots,0,1,0,\ldots,0)$ denotes the~$k$-th standard unit vector with 1 on the~$k$-th entry.
These facts are already known; cf.~\cite[Theorem 23.2]{Ismail} and~\cite{HV}.
Moreover, if $r=1$, the type II multiple orthogonal polynomials reduce to the usual orthogonal polynomials, hence,
Theorems~\ref{thm:l even} and~\ref{thm:l odd} generalize classical results of Karlin and Szeg\H{o} mentioned at the
beginning.

Finally, one can also consider the Wronskians involving type I multiple orthogonal polyno\-mials by replacing
$P_{\vec{n}}$ in~\eqref{def:wronskian} by $Q_{\vec{n}}$ def\/ined in~\eqref{def:Qn}.
In this case, the Wronskians are not polynomials in general.

\begin{Theorem}
\label{thm:type 1}
Let $w_1,w_2,\ldots,w_r$ be~$r$ sufficiently many times differentiable weights as in Theo\-rem~{\rm \ref{thm:l even}} and assume
that~$l$ is even.
For each path starting from the multi-index $\vec{n}$, we have that the function
$W(Q_{\vec{n}}(x),Q_{\vec{n}_1}(x),\ldots,Q_{\vec{n}_{l-1}}(x))$ given in~\eqref{def:wronskian} keeps
a~constant sign on $(a,b)$, if the Wronskian is well defined.
\end{Theorem}

\subsection{Outline of the paper}
The rest of this paper is organized as follows.
In Section~\ref{sec:proof}, we prove Theorems~\ref{thm:l even}--\ref{thm:type 1}.
We next give two applications of our main results.
In Section~\ref{sec:moments}, we show that the~$m$-th moment of the average characteristic polynomials for multiple
orthogonal polynomial ensembles can be expressed as a~Wronskian of the type II multiple orthogonal polynomials.
It is then straightforward to obtain the distinct behavior of the moments for odd and even~$m$ in a~special multiple
orthogonal ensemble~-- the AT ensemble.
In Section~\ref{sec:inequality} we derive the inequalities of Tur\'an type for some classical multiple orthogonal
polynomials, namely, for multiple Hermite polynomials and multiple Laguerre polynomials.
We conclude this paper with numerical study of the zero conf\/igurations of the Wronskians for multiple Hermite
polynomials and multiple Laguerre polynomials in Section~\ref{sec:zeros}.
It comes out that the zeros have fascinating   and regular conf\/igurations in the complex plane, which might be of
independent interest.

\section{Proofs of Theorems~\ref{thm:l even}--\ref{thm:type 1}}
\label{sec:proof}

We shall prove our main theorems by extending the arguments in~\cite{KS}.
Roughly speaking, the proofs of Theorems~\ref{thm:l even} and~\ref{thm:type 1} use the properties of an AT system, while
for the proof of Theorem~\ref{thm:l odd} one needs Theorem~\ref{thm:l even} and Sylvester's theorem concerning the
relation between the determinants of a~square matrix and its minors.

\subsection{Proof of Theorem~\ref{thm:l even}}

We f\/irst show that the Wronskian~\eqref{eq: l even} keeps a~constant sign on the real axis for~$l$ even.
If this is not true, we may f\/ind a~real number $x_0$ such that $W(\vec{n},l;x_0)=0$.
This in turn implies the existence of the constants $\lambda_0,\lambda_1,\ldots, \lambda_{l-1}$ depending on the path
such that the function
\begin{gather}
\label{def:f}
f(x):=\sum\limits_{i=0}^{l-1}\lambda_iP_{\vec{n}_i}(x)
\end{gather}
satisf\/ies
\begin{gather*}
f^{(k)}(x_0)=0,
\qquad
k=0,1,\ldots,l-1.
\end{gather*}
Thus, $x_0$ is a~zero of $f(x)$ of multiplicity at least~$l$.

We further claim that~$f$ has at least $|\vec{n}|$ zeros on $(a,b)$ where~$f$ changes sign.
Such zeros are also called {\it nodal zeros}.
To see this, we f\/irst observe from~\eqref{eq:orthogonality of Pn} and~\eqref{def:f} that the equality
\begin{gather}
\label{eq:orth of f}
\int f(x)\sum\limits_{i=1}^{r}q_{i}(x)w_{i}(x)\,dx=0
\end{gather}
holds for any polynomials $q_i(x)$ of degree less than or equal to $n_i-1$, where we also make use of the fact that $
\vec{n}\leq \vec{n}_1\leq \dots \leq \vec{n}_{l-2}\leq\vec{n}_{l-1}$.
In order to obtain a~contradiction, suppose that~$f$ has at most $k\leq |\vec{n}|-1$ nodal zeros on $(a,b)$, say,
$x_1,x_2, \ldots,x_k$, then
\begin{gather}
\label{eq:factor f}
f(x)=\prod\limits_{i=1}^{k}(x-x_i)Q(x),
\end{gather}
where~$Q$ does not change the sign on the the interval $(a,b)$.
Since $k\leq |\vec{n}|-1$, it is always possible to construct a~multi-index $\vec{m}=(m_1,m_2,\ldots,
m_r)\in\mathbb{N}^r$ such that $\vec{m} \leq \vec{n}$, and $|\vec{m}|=k+1 \leq |\vec{n}|$ for any given initial
multi-index $\vec{n}$.
Then there exist polynomials $\widetilde q_i$, $1\leq i\leq r$, with degrees less than or equal to $m_i-1$ satisfying
\begin{gather*}
\sum\limits_{i=1}^{r}\widetilde q_{i}(x)w_{i}(x)=
\begin{cases}
0, & \text{if} \quad x=x_1,x_2,\ldots,x_k,
\\
1, & \text{if} \quad x=x_{k+1},
\end{cases}
\end{gather*}
where $x_1,x_2,\ldots,x_k$ are the nodal zeros as in~\eqref{eq:factor f} and $x_{k+1}$ is an arbitrary point on $(a,b)$
but dif\/ferent from those~$k$ nodal zeros.
Indeed, this is equivalent to solving a~linear system of $|\vec{m}|=k+1$ equations for the unknown coef\/f\/icients of
$\widetilde q_i$.
This system is uniquely solvable if the matrix
\begin{gather*}
\left(
\begin{matrix}
w_1(x_1) & \cdots & x_1^{m_1-1}w_1(x_1)
\\
w_1(x_2) & \cdots & x_2^{m_1-1}w_1(x_2)
\\
\vdots & \vdots & \vdots
\\
w_1(x_{k+1}) & \cdots & x_{k+1}^{m_1-1}w_1(x_{k+1})
\end{matrix}
\cdots \cdots
\begin{matrix}
w_r(x_1) & \cdots & x_1^{m_r-1}w_r(x_1)
\\
w_r(x_2) & \cdots & x_2^{m_r-1}w_r(x_2)
\\
\vdots & \vdots & \vdots
\\
w_r(x_{k+1}) & \cdots & x_{k+1}^{m_r-1}w_r(x_{k+1})
\end{matrix}
\right)
\end{gather*}
is not singular, which is immediate on account of the fact that our weights $(w_1,w_2,\ldots,w_r)$ form an AT system on
$[a,b]$ for all the multi-indices.
The Chebyshev property further indicates that the function
$\prod\limits_{i=1}^{k}(x-x_i)(\sum\limits_{i=1}^{r}\widetilde q_{i}(x)w_{i}(x))$ does not change the sign on
$(a,b)$ (cf.~\cite{Ker}).
This, together with~\eqref{eq:factor f}, implies that
\begin{gather*}
\int f(x)\sum\limits_{i=1}^{r}\widetilde q_{i}(x)w_{i}(x)\,dx\neq 0,
\end{gather*}
which is a~contradiction with~\eqref{eq:orth of f}.
Thus, we have proved~$f$ has at least $|\vec{n}|$ nodal zeros on $(a,b)$.

Note that our def\/inition of~$f$ in~\eqref{def:f} shows that~$f$ is a~polynomial of degree less than or equal to
$|\vec{n}|+l-1$.
If $x_0$ is dif\/ferent from all these nodal zeros, then~$f$ will have at least $|\vec{n}|+l$ zeros, which is
a~contradiction.
On the other hand, if $x_0$ coincides with one of the nodal zeros, multiplicity of $x_0$ must be at least $l+1$ (since
the zero of even multiplicity~$l$ cannot be nodal), hence, the number of total zeros is at least
$|\vec{n}|-1+l+1=|\vec{n}|+l$, again a~contradiction.
In summary, we have proved that $W(\vec{n},l;x)$ keeps a~constant sign on the real axis if~$l$ is even.

{\sloppy Finally, the inequality~\eqref{eq: l even} follows from the observation that the leading term of the $W(\vec{n},l;x)$ is $x^{l|\vec{n}|}$,
up to a~positive constant.
This can be seen by f\/irst multiplying the rows of~\eqref{def:wronskian} by $1,x,\ldots,x^{l-1}$ and the columns by
$x^{-|\vec{n}|},x^{-|\vec{n}|-1},\ldots, x^{-|\vec{n}|-l+1}$, respectively, and then letting $x\to\infty$.

}

\subsection{Proof of Theorem~\ref{thm:l odd}
} The proof of this theorem is similar to that for the case of orthogonal polynomials in~\cite{KS}.
We start with a~special form of Sylvester's identity which states that for any square matrix~$A$ of size~$n$ and $n \geq
m_1 \geq m_2 \geq 1$ and $n \geq n_1>n_2 \geq 1$, the following identity (cf.~\cite{Gan}) holds:
\begin{gather}
\label{eq:sylvester identity}
\det A\cdot \det A[m_1,m_2;n_1,n_2]= \det
\begin{pmatrix}
\det A[m_1;n_1] & \det A[m_1;n_2]
\\
\det A[m_2;n_1] & \det A[m_2;n_2]
\end{pmatrix}
,
\end{gather}
where $A[m_1,m_2;n_1,n_2]$ denotes the submatrix obtained from~$A$ by deleting rows $m_1$, $m_2$ and columns $n_1$, $n_2$, and
a~similar def\/inition holds for $A[m_i;n_j]$, $i,j=1,2$.

Given any path, say, $(\vec{n}_0,\vec{n}_1,\ldots,\vec{n}_{l-1})$, we can f\/ind another path
$(\vec{n}_1,\vec{n}_2,\ldots,\vec{n}_{l})$ whose f\/irst $l-1$ multi-indices coincide with the last $l-1$ multi-indices of
the original path.
Applying~\eqref{eq:sylvester identity} to $\det A= W(\vec{n}_0,\vec{n}_1,\ldots,\vec{n}_{l-1},\vec{n}_{l};x)$ with
$m_1=l+1$, $m_2=l$ and $n_1=l+1$, $n_2=1$, we have
\begin{gather*}
 W(\vec{n}_0,\ldots,\vec{n}_{l};x) \cdot W(\vec{n}_1,\ldots,\vec{n}_{l-1};x)
 =\det
\begin{pmatrix}
W(\vec{n}_0,\ldots,\vec{n}_{l-1};x) & W(\vec{n}_1,\ldots,\vec{n}_{l};x)
\\
W'(\vec{n}_0,\ldots,\vec{n}_{l-1};x) & W'(\vec{n}_1,\ldots,\vec{n}_{l};x)
\end{pmatrix},
\end{gather*}
where the derivative $'$ is with respect to~$x$.
Recall that~$l$ is odd, hence, Theorem~\ref{thm:l even} gives us
\begin{gather*}
W(\vec{n}_0,\ldots,\vec{n}_{l-1};x) W'(\vec{n}_1,\ldots,\vec{n}_{l};x)
-W'(\vec{n}_0,\ldots,\vec{n}_{l-1};x)W(\vec{n}_1,\ldots,\vec{n}_{l};x) >0,
\end{gather*}
for all $x\in\mathbb{R}$.
This inequality shows that all zeros of $W(\vec{n}_0,\ldots,\vec{n}_{l-1};x)$ must be simple.
Moreover, we have
\begin{gather}
\label{eq:inequality 1}
W(\vec{n}_0,\ldots,\vec{n}_{l-1};x) W'(\vec{n}_1,\ldots,\vec{n}_{l};x)>0,
\qquad
\textrm{if}
\qquad
W(\vec{n}_1,\ldots,\vec{n}_{l};x)=0,
\end{gather}
and
\begin{gather}
\label{eq:inequality 2}
W'(\vec{n}_0,\ldots,\vec{n}_{l-1};y)W(\vec{n}_1,\ldots,\vec{n}_{l};y)<0,
\qquad
\textrm{if}
\qquad
W(\vec{n}_0,\ldots,\vec{n}_{l-1};y)=0.
\end{gather}
If $x_0$ and $\widetilde x_0$ are two consecutive zeros of $W(\vec{n}_1,\ldots,\vec{n}_{l};x)$, then
\begin{gather*}
W'(\vec{n}_1,\ldots,\vec{n}_{l};x_0)W'(\vec{n}_1,\ldots,\vec{n}_{l};\widetilde x_0)<0.
\end{gather*}
This, together with~\eqref{eq:inequality 1}, implies that
\begin{gather*}
W(\vec{n}_0,\ldots,\vec{n}_{l-1};x_0)W(\vec{n}_0,\ldots,\vec{n}_{l-1};\widetilde x_0)<0.
\end{gather*}
Hence, there must be a~zero of $W(\vec{n}_0,\ldots,\vec{n}_{l-1};x)$ between $(x_0,\widetilde x_0)$.
By~\eqref{eq:inequality 2} and the same argument, we conclude that there exists at least one zero of
$W(\vec{n}_1,\ldots,\vec{n}_{l};x)$ between any two consecutive zeros of $W(\vec{n}_0,\ldots,\vec{n}_{l-1};x)$.
This completes the proof of the simplicity and the interlacing property of real zeros stated in Theorem~\ref{thm:l odd}.

Finally, it remains to calculate the number of real zeros of the Wronskian.
This can be achieved by induction argument on $|\vec{n}|$.
If $|\vec{n}|=0$ (i.e., $\vec{n}=\vec{0}$), then $W(\vec{0},l;x)>0$, since the Wronskian matrix reduces to the upper
diagonal matrix with positive diagonal entries.
Hence, there is no real zero in this case.
If $|\vec{n}|=1$, we have $W(\vec{n},l;x)\sim x^{|\vec{n}|l}=x^l$ as $x\to \pm \infty$ and~$l$ is odd, thus,
$W(\vec{n},l;x)$ has at least one real zero.
If it has more than one real zero, by interlacing property this will lead to the existence of a~real zero for certain
Wronskian associated with a~path starting from $\vec{0}$, hence, a~contradiction.
Suppose now that $W(\vec{n},l;x)$ has exactly $k>0$ simple real zeros if $|\vec{n}|=k$.
When $|\vec{n}|=k+1$, for any path starting from $\vec{n}$, we can f\/ind another path starting from $\vec{m}$ such that
$|\vec{m}|=k$ and the zeros of $W(\vec{n},l;x)$ and $W(\vec{m},l;x)$ interlace.
Then, $W(\vec{n},l;x)$ will have at least $k-1$ simple zeros.
If $x_k$ is the largest zero of $W(\vec{m},l;x)$, then $W'(\vec{m},l;x_k)>0$ since $W(\vec{m},l;x)$ is positive for~$x$
large.
From~\eqref{eq:inequality 2}, we have that $W(\vec{n},l;x_k)<0$, hence there will be at least one zero of
$W(\vec{n},l;x)$ on the right hand side of~$x_k$.
There would be only one such zero, again by interlacing property.
Similar argument implies that there will be exactly one zero on the left of the smallest real zero of $W(\vec{m},l;x)$.
Thus, $W(\vec{n},l;x)$ will have exactly $k+1=|\vec{n}|$ real simple zeros.

\subsection{Proof of Theorem~\ref{thm:type 1}
} The proof is similar to that of Theorem~\ref{thm:l even}.
Suppose that there exists a~point $x_0\in(a,b)$ such that
$W(Q_{\vec{n}}(x_0),Q_{\vec{n}_1}(x_0),\ldots,Q_{\vec{n}_{l-1}}(x_0))=0$.
Then we can f\/ind constants $\lambda_i$, $i=0,1,\ldots,l-1$, such that the function
\begin{gather*}
f(x):=\sum\limits_{i=0}^{r-1}\lambda_i Q_{\vec{n}_i}(x),
\qquad
\vec{n}_0=\vec{n},
\end{gather*}
has a~zero at $x_0$ of multiplicity at least~$l$.
Since
\begin{gather*}
\int f(x)p(x)\,dx=0
\end{gather*}
for any polynomial~$p$ of degree less than or equal to $|\vec{n}|-2$, we have that~$f$ has at least $|\vec{n}|-1$ nodal
zeros.
Thus, as in the proof of Theorem~\ref{thm:l even}, we conclude that~$f$ will have at least $|\vec{n}|+l-1$ zeros.
This is a~contradiction to the fact that~$f$ is a~linear combination of the Chebyshev system for the multi-index
$\vec{n}_{l-1}$, which has at most $|\vec{n}_{l-1}|-1=|\vec{n}|+l-2$ zeros on $[a,b]$.

\section{Moments of the average characteristic polynomials\\ for multiple orthogonal polynomial ensembles}
\label{sec:moments}

In this section, we shall apply our results to the moments of the average characteristic polynomials for multiple
orthogonal polynomial ensembles.

It is well-known that, besides the interest from the approximation theory, multiple orthogonal polynomials have also
arisen recently in a~natural way in certain models of mathematical physics, including random matrix theory,
non-intersecting paths, etc; cf.~\cite{Kui1,Kui2} and the re\-fe\-ren\-ces therein.
Indeed, let us consider $|\vec{n}| = \sum\limits_{i=1}^{r}n_i= n$ random points on the real line whose joint probability
density function (p.d.f.) can be written as a~product of two determinants:
\begin{gather}
\label{def:pjdf}
\frac{1}{Z_n}\det\left(f_i(x_j)\right)_{i,j=1}^n\det\left(g_i(x_j)\right)_{i,j=1}^n,
\end{gather}
where $Z_n$ is the normalizing constant to make the total probability on $\mathbb{R}^n$ equal to one, and the two
sequences of functions $f_i$, $g_i$ are given~by
\begin{gather*}
f_i(x)=x^{i-1},
\qquad
i=1,\ldots, n,
\end{gather*}
and
\begin{gather*}
  g_{i}(x)= x^{i-1}w_1(x),
\qquad
i=1,\ldots,n_1,
\\
  g_{n_1+i}(x)=x^{i-1}w_2(x),
\qquad
i=1,\ldots,n_2,
\\
 \cdots\cdots\cdots\cdots\cdots\cdots\cdots\cdots\cdots\cdots\cdots\cdots\cdots\cdots\cdots\cdots
\\
  g_{n_1+\dots+n_{r-1}+i}(x)=x^{i-1}w_r(x),
\qquad
i=1,\ldots,n_r.
\end{gather*}
We call this stochastic model a~multiple orthogonal polynomial ensemble, since one has
\begin{gather}
P_{\vec{n}}(z) =\mathbb{E}\left[\prod\limits_{k=1}^n(z-x_k)\right],
\qquad
z\in\mathbb{C},
\label{eq:acp}
\end{gather}
where the expectation $\mathbb{E}$ is taken with respect to the p.d.f.~\eqref{def:pjdf}.

The formula~\eqref{eq:acp} tells us that the type II multiple orthogonal polynomial $P_{\vec{n}}(z)$ can be viewed as
the average of the random polynomials $\prod\limits_{k=1}^n(z-x_k)$ whose roots are distributed according
to~\eqref{def:pjdf}.
As a~consequence, $P_{\vec{n}}$ is also called the average characteristic polynomial if the
distribution~\eqref{def:pjdf} can be interpreted as the particle distribution of certain stochastic models.
In particular, we mention that one example falling into this category is from the random matrix model with external
source~\cite{BH981,BH982,Zinn}, which was f\/irst observed in~\cite{BleKuij2}.

The result~\eqref{eq:acp} is further extended by Delvaux in~\cite{Delvaux} to arbitrary products and ratios of
characteristic polynomials, which are def\/ined~by
\begin{gather*}
\mathbb{E}\left[\frac{\prod\limits_{j=1}^m\prod\limits_{k=1}^n(z_j-x_k)}{\prod\limits_{j=1}^l\prod\limits_{k=1}^n(y_j-x_k)}\right],
\qquad
l,m, n\in\mathbb{N},
\end{gather*}
where $z_1,\ldots,z_m \in \mathbb{C}$, $y_1,\ldots, y_l \in \mathbb{C} \setminus \mathbb{R}$, and all the numbers in the
set $(z_1,\ldots,z_m,
y_1,\ldots, y_l)$ are pairwise dif\/ferent.
It turns out that the products/ratios of the average characteristic polynomials for multiple orthogonal polynomial
ensembles can be expressed as the determinants whose entries involve the blocks of the Riemann--Hilbert matrix
characterizing multiple orthogonal polynomials and a~matrix-valued version of the Christof\/fel--Darboux kernel.
Particularly, in case $l=0$, it follows from~\cite[Theorem 1.8]{Delvaux} that
\begin{gather}
\mathbb{E}\left[\prod\limits_{j=1}^m\prod\limits_{k=1}^n\!(z_j-x_k)\right]=\frac{1}{\prod\limits_{1\leq i<j \leq m}(z_j-z_i)}
\det \!\left(
\begin{matrix}
P_{\vec{n}_0}(z_1) & P_{\vec{n}_{1}}(z_1) & \cdots & P_{\vec{n}_{m-1}}(z_1)
\\
P_{\vec{n}_0}(z_2) & P_{\vec{n}_{1}}(z_2) & \cdots & P_{\vec{n}_{m-1}}(z_2)
\\
\vdots & \vdots & \vdots & \vdots
\\
P_{\vec{n}_0}(z_m) & P_{\vec{n}_{1}}(z_m) & \cdots & P_{\vec{n}_{m-1}}(z_m)
\end{matrix}
\right)\!,\!\!\!
\label{eq:prod}
\end{gather}
where $\vec{n}_0=\vec{n}$, $\vec{n}_k\geq \vec{n}_{k-1}$ componentwise and $|\vec{n}_k-\vec{n}_{k-1}|=1$ for
$k=1,\ldots,m$, which is actually an arbitrary path connecting $\vec{n}_0$ to $\vec{n}_{m-1}$; see the def\/inition at the
beginning of Section~\ref{subsec:main results}.

Now, let all $z_j$ in~\eqref{eq:prod} tend to~$z$, an algebraic manipulation (cf.~\cite[Theorem 1.2.4]{Hua}) shows that
the moments of the average characteristic polynomials for multiple orthogonal polynomial ensembles can be expressed as
the Wronskians of type II multiple orthogonal polynomials:
\begin{gather}
\label{eq:moments}
\mathbb{E}\left[\prod\limits_{k=1}^n(z-x_k)^m\right]=\frac{1}{\prod\limits_{i=0}^{m-1}i!} W\left(P_{\vec{n}_0}(z),
P_{\vec{n}_{1}}(z), \ldots, P_{\vec{n}_{m-1}}(z)\right),
\end{gather}
where~$W$ is def\/ined in~\eqref{def:wronskian}.
Note that when $g_i(x)=x^{i-1}w(x)$ in~\eqref{def:pjdf} (i.e., in the case of orthogonal polynomial ensembles), the
formula above was f\/irst shown by Br\'{e}zin and Hikami~\cite{BH2001}; see also~\cite{MN}.

Combining~\eqref{eq:moments} and Theorems~\ref{thm:l even}--\ref{thm:l odd}, we obtain
\begin{Corollary}
Assume that the weights $(w_1,w_2,\ldots,w_r)$ in~\eqref{def:pjdf} form an AT system on $[a,b]$ for all the
multi-indices in $\mathbb{N}^r$ $($i.e., an AT ensemble in the sense of~{\rm \cite{Kui1})}.
Then we have that the moments of the average characteristic polynomials with respect to~\eqref{def:pjdf}
\begin{gather*}
\mathbb{E}\left[\prod\limits_{k=1}^n(z-x_k)^m\right]
\end{gather*}
are strictly positive on the real axis if~$m$ is even; while for odd~$m$, the moments admit oscillatory behavior as
stated in Theorem~{\rm \ref{thm:l odd}}.
\end{Corollary}

\section{Some inequalities for multiple orthogonal polynomials}
\label{sec:inequality}

In this section, we shall use our results to derive the inequalities of Tur\'an type for some classical multiple
orthogonal polynomials, namely, for multiple Hermite polynomials and multiple Laguerre polynomials.
It is known that the weights for these polynomials form an AT system for any multi-index $\vec{n}\in\mathbb{N}^r$.

\subsection{Tur\'an inequalities for multiple Hermite polynomials}
Multiple Hermite polynomials are def\/ined~by
\begin{gather*}
\int_{-\infty}^\infty x^k H_{\vec{n}}(x) e^{-x^2+c_jx}\,dx =0,
\qquad
k = 0, 1, \ldots, n_j-1,
\end{gather*}
for $j=1,\ldots,r$, where $c_i \neq c_j$ if $i \neq j$; cf.~\cite{BleKuij},~\cite[\S~23.5]{Ismail} and~\cite[\S~3.4]{WVAEC}.
An explicit formula for multiple Hermite polynomials is
\begin{gather*}
H_{\vec{n}}=\frac{(-1)^{|\vec{n}|}}{2^{|\vec{n}|}}\sum\limits_{k_1=0}^{n_1}\cdots
\sum\limits_{k_r=0}^{n_r}\binom{n_1}{k_1}\cdots\binom{n_r}{k_r}
\prod\limits_{j=1}^{r}c_j^{n_j-k_j}(-1)^{|\vec{k}|}H_{|\vec{k}|}(x),
\end{gather*}
where $\vec{k}=(k_1,\ldots,k_r)$ and $H_{|\vec{k}|}$ is the usual Hermite polynomial of degree $|\vec{k}|$ with the leading
coef\/f\/icient $2^{|\vec{k}|}$.
The following statement holds for multiple Hermite polynomials.
\begin{Theorem}
The multiple Hermite polynomials satisfy the following inequalities:
\begin{gather}
\label{eq:turan MHP}
H_{\vec{n}+\vec{e}_j}(x)H_{\vec{n}+\vec{e}_k}(x)- H_{\vec{n}}(x)H_{\vec{n}+\vec{e}_j+\vec{e}_k}(x)>0,
\qquad
x\in\mathbb{R},
\end{gather}
for $j,k=1,\ldots,r$, where $\vec{e}_i = (0,\ldots,0,1,0,\ldots,0)$ denotes the~$i$-th standard unit vector with~$1$ on
the~$i$-th entry.
In particular, by taking $j=k$, we have
\begin{gather}
\label{eq:j=k MHP}
H_{\vec{n}+\vec{e}_j}^2(x)- H_{\vec{n}}(x)H_{\vec{n}+2\vec{e}_j}(x)>0,
\qquad
x\in\mathbb{R}.
\end{gather}
\end{Theorem}

\begin{Remark}
It is readily seen that the formula~\eqref{eq:j=k MHP} extends the Tur\'an inequality for the Hermite
polynomials~\cite{KS}.
\end{Remark}

\begin{proof}
From Theorem~\ref{thm:l even} with $l=2$, it follows that
\begin{gather}
\label{eq:MHP inequality}
\det
\begin{pmatrix}
H_{\vec{n}}(x) & H_{\vec{n}+\vec{e}_j}(x)
\\
H_{\vec{n}}'(x) & H_{\vec{n}+\vec{e}_j}'(x)
\end{pmatrix}
>0,
\qquad
j=1,\ldots,r,
\end{gather}
for any $x\in\mathbb{R}$.
Note that $H_{\vec{n}}$ satisf\/ies the following raising operations (cf.~\cite[\S~23.8.2]{Ismail})
\begin{gather*}
\frac{d}{dx}\big(e^{-x^2+c_k x}H_{\vec{n}}(x)\big)=-2e^{-x^2+c_k x}H_{\vec{n}+\vec{e}_k}(x),
\qquad
k=1,\ldots,r,
\end{gather*}
or, equivalently,
\begin{gather}
\label{eq:rasing operators}
\frac{d}{dx}H_{\vec{n}}(x)=-2H_{\vec{n}+\vec{e}_k}(x)+(2x-c_k)H_{\vec{n}}(x),
\qquad
k=1,\ldots,r.
\end{gather}
Inserting this formula into~\eqref{eq:MHP inequality} gives~\eqref{eq:turan MHP}.
\end{proof}

It is worthwhile to point out that the inequality~\eqref{eq:turan MHP} is independent of the parameters $c_j$ appearing
in the weight functions.
Furthermore, by choosing the path in the Wronskian matrix to be $(\vec{n}, \vec{n}+\vec{e}_j, \ldots,
\vec{n}+(l-1)\vec{e}_j)$ for f\/ixed $j=1,\ldots,r$ and using~\eqref{eq:rasing operators} with $k=j$, it is readily seen
that
\begin{gather}
 W\big(H_{\vec{n}}(x), H_{\vec{n}+\vec{e}_j}(x), \ldots, H_{\vec{n}+(l-1)\vec{e}_j}(x)\big)
\nonumber
\\
\qquad
 = (-2)^{\frac{l(l-1)}{2}} \det \left(
\begin{matrix}
H_{\vec{n}}(x) & H_{\vec{n}+\vec{e}_j}(x) & \dots & H_{\vec{n}+(l-1)\vec{e}_j}(x)
\\
H_{\vec{n}+\vec{e}_j}(x) & H_{\vec{n}+2\vec{e}_j}(x) & \dots & H_{\vec{n}+l\vec{e}_j}(x)
\\
\vdots & \vdots & \vdots & \vdots
\\
H_{\vec{n}+(l-1)\vec{e}_j}(x) & H_{\vec{n}+l\vec{e}_j}(x) & \dots & H_{\vec{n}+2(l-1)\vec{e}_j}(x)
\end{matrix}
\right),
\label{eq:wroskian to Hankel}
\end{gather}
that is, we pass from the determinant of the Wronskian type to the Hankel determinant.
This, together with Theorem~\ref{thm:l even}, implies
\begin{Corollary}
Let $T(H_{\vec{n}}(x), H_{\vec{n}+\vec{e}_j}(x), \ldots, H_{\vec{n}+(l-1)\vec{e}_j}(x))$ be the Hankel
determinant of multiple Hermite polynomials on the right hand side of~\eqref{eq:wroskian to Hankel}, then
\begin{gather*}
(-1)^{\frac{l(l-1)}{2}}T\big(H_{\vec{n}}(x), H_{\vec{n}+\vec{e}_j}(x), \ldots, H_{\vec{n}+(l-1)\vec{e}_j}(x)\big)>0,
\qquad
x\in\mathbb{R},
\qquad
j=1,\ldots,r,
\end{gather*}
if~$l$ is even.
\end{Corollary}

\subsection{Two-parameter Tur\'an inequalities for multiple Laguerre polynomials}

There are two kinds of multiple Laguerre polynomials.
Multiple Laguerre polynomials of the f\/irst kind are def\/ined by the orthogonality conditions
\begin{gather*}
\int_0^\infty x^k L_{\vec{n}}^{\vec{\alpha}}(x) x^{\alpha_j} e^{-x}\,dx =0,
\qquad
k = 0, 1, \ldots, n_j-1,
\end{gather*}
for $j=1,\ldots,r$, where $\vec{\alpha}=(\alpha_1,\ldots,\alpha_r)$ with $\alpha_j>-1$ and $\alpha_i-\alpha_j\notin
\mathbb{Z}$ whenever $i \neq j$; cf.~\cite{BleKuij},~\cite[\S~23.4.1]{Ismail} and~\cite[\S~3.2]{WVAEC}.
An explicit formula for the multiple Laguerre polynomials of the f\/irst kind is
\begin{gather}
L_{\vec{n}}^{\vec{\alpha}}(x)=\sum\limits_{k_1=0}^{n_1} \cdots \sum\limits_{k_r=0}^{n_r} \binom{n_1}{k_1} \cdots
\binom{n_r}{k_r} \binom{n_r+\alpha_r}{k_r}\cdots\nonumber
\\
\hphantom{L_{\vec{n}}^{\vec{\alpha}}(x)=}{}
\times \binom{|\vec{n}|-|\vec{k}|+k_1+\alpha_1}{k_1} \prod\limits_{i=1}^{r}k_i! (-1)^{|\vec{k}|}x^{|\vec{n}|-|\vec{k}|}.\label{eq:ML1st}
\end{gather}

Multiple Laguerre polynomials of the second kind are def\/ined by the orthogonality conditions
\begin{gather*}
\int_0^\infty x^k L_{\vec{n}}^{(\alpha,\vec{c})}(x) x^{\alpha} e^{-c_j x}\,dx =0,
\qquad
k = 0, 1, \ldots, n_j-1,
\end{gather*}
for $j=1,\ldots,r$, where $\vec{c}=(c_1,\ldots,c_r)$ and we assume that $\alpha>-1$, $c_j>0$ and $c_i \neq c_j$ whenever $i \neq j$;
cf.~\cite{BleKuij},
\cite[\S~23.4.2]{Ismail},
\cite[Remark~5 on p.~160]{NikSor} and~\cite[\S~3.3]{WVAEC}.
An explicit formula for these polynomials is
\begin{gather*}
L_{\vec{n}}^{(\alpha,\vec{c})}(x)=\sum\limits_{k_1=0}^{n_1} \cdots \sum\limits_{k_r=0}^{n_r} \binom{n_1}{k_1} \cdots
\binom{n_r}{k_r} \binom{|\vec{n}|+\alpha}{|\vec{k}|}|\vec{k}|!\frac{(-1)^{|\vec{k}|}}
{\prod\limits_{j=1}^{r}c_j^{k_j}}x^{|\vec{n}|-|\vec{k}|}.
\end{gather*}

\begin{figure}[t]
\centering
\includegraphics[width=7.5cm]{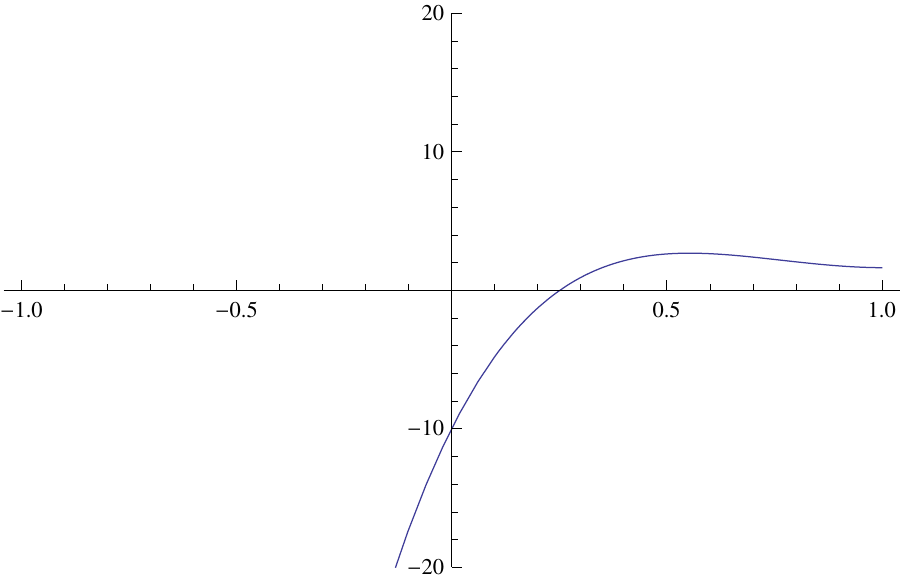}

\caption{Plot of~\eqref{eq:couterex} with $\vec{n}=(n,m)=(1,1)$ and
$\vec{\alpha}=(1/2,1/3)$.}
\label{fig:0}
\end{figure}

It turns out that the Tur\'an inequalities of the form~\eqref{eq:turan MHP} do not hold for multiple Laguerre
polynomials in general.
We can easily calculate from~\eqref{eq:ML1st} that, for instance, for multiple Laguerre polynomials of the f\/irst kind
with $r=2$, $\vec{n}=(n,m)=(1,1)$ and $\vec{\alpha}=(1/2,1/3)$ the expression
\begin{gather}
\label{eq:couterex}
\big(L^{\vec{\alpha}}_{n+1,m}(x)\big)^2-L^{\vec{\alpha}}_{n,m}(x)L^{\vec{\alpha}}_{n+2,m}(x)
\end{gather}
is reduced to
\begin{gather*}
2x^5-\frac{119}{6}x^4+\frac{647}{9}x^3-\frac{7495}{72}x^2+\frac{185}{3}x-10,
\end{gather*}
which can be both positive and negative; see Fig.~\ref{fig:0}.
We can also f\/ind similar counterexamples for the other choices of indices and also for multiple Laguerre polynomials of
the second kind.
However, multiple Laguerre polynomials satisfy the following two-parameter Tur\'an inequalities in the sense
of~\cite{Busto}.

\begin{Theorem}
\label{thm:turan MLP}
For multiple Laguerre polynomials of the first kind we have
\begin{gather}
\label{eq:turan MLP1}
L_{\vec{n}+\vec{e}_k}^{\vec{\alpha}}(x)L_{\vec{n}+\vec{e}_j}^{\vec{\alpha}-\vec{e}_j}(x)-
L_{\vec{n}}^{\vec{\alpha}}(x)L_{\vec{n}+\vec{e}_j+\vec{e}_k}^{\vec{\alpha}-\vec{e}_j}(x)>0,
\qquad
x>0,
\end{gather}
for $\vec{\alpha}>\vec{0}$ and $j,k=1,\ldots,r$.

Similarly, for multiple Laguerre polynomials of the second kind, we have
\begin{gather}
\label{eq:turan MLP2}
L_{\vec{n}+\vec{e}_k}^{(\alpha,\vec{c})}(x)L_{\vec{n}+\vec{e}_j}^{(\alpha-1,\vec{c})}(x)-
L_{\vec{n}}^{(\alpha,\vec{c})}(x)L_{\vec{n}+\vec{e}_j+\vec{e}_k}^{(\alpha-1,\vec{c})}(x)>0,
\qquad
x>0,
\end{gather}
for $\alpha>0$ and $j,k=1,\ldots,r$.
\end{Theorem}
\begin{proof}
From Theorem~\ref{thm:l even} with $l=2$, we see that
\begin{gather}
\label{eq:MLP inequality}
\det
\begin{pmatrix}
L_{\vec{n}}(x) & L_{\vec{n}+\vec{e}_k}(x)
\\
xL_{\vec{n}}'(x) & xL_{\vec{n}+\vec{e}_k}'(x)
\end{pmatrix}
>0,
\qquad
x>0,
\end{gather}
for $k=1,\ldots,r$, where $L_{\vec{n}}$ is the multiple Laguerre polynomial of the f\/irst or second kind.
Note that $L_{\vec{n}}^{\vec{\alpha}}$ satisf\/ies the following raising operator (cf.~\cite[\S~23.4.2]{Ismail}):
\begin{gather*}
\frac{d}{dx}\big(x^{\alpha_j}e^{-x}L_{\vec{n}}^{\vec{\alpha}}(x)\big)=-x^{\alpha_j-1}e^{-x}
L_{\vec{n}+\vec{e}_j}^{\vec{\alpha}-\vec{e}_j}(x),
\qquad
j=1,\ldots,r,
\end{gather*}
or, equivalently,
\begin{gather*}
x\frac{d}{dx}L_{\vec{n}}^{\vec{\alpha}}(x)=(x-\alpha_j)L_{\vec{n}}^{\vec{\alpha}}(x)
-L_{\vec{n}+\vec{e}_j}^{\vec{\alpha}-\vec{e}_j}(x),
\qquad
j=1,\ldots,r,
\end{gather*}
Substituting this formula into~\eqref{eq:MLP inequality} gives us~\eqref{eq:turan MLP1}.

Similarly,~\eqref{eq:turan MLP2} follows from~\eqref{eq:MLP inequality} and the relation (cf.~\cite[\S~23.4.4]{Ismail})
\begin{gather*}
\frac{d}{dx}\big(x^{\alpha}e^{-c_j x}L_{\vec{n}}^{(\alpha,\vec{c})}(x)\big)=-c_jx^{\alpha-1}e^{-c_jx}
L_{\vec{n}+\vec{e}_j}^{(\alpha-1,\vec{c})}(x),
\qquad
j=1,\ldots,r,
\end{gather*}
or, equivalently,
\begin{gather*}
x\frac{d}{dx}L_{\vec{n}}^{(\alpha,\vec{c})}(x)=(c_jx-\alpha)L_{\vec{n}}^{(\alpha,\vec{c})}(x)
-c_jL_{\vec{n}+\vec{e}_j}^{(\alpha-1,\vec{c})}(x),
\qquad
j=1,\ldots,r.
\end{gather*}
This completes the proof of Theorem~\ref{thm:turan MLP}.
\end{proof}

\section{Conf\/igurations of zeros for the Wronskians\\ of multiple Hermite and Laguerre polynomials}
\label{sec:zeros}

We conclude this paper by studying numerically the geometric conf\/iguration of zeros of the Wronskians of multiple
Hermite and Laguerre (of both kinds) polynomials using {\it Mathematica}\footnote{\url{http://www.wolfram.com}}.
Our motivation arises from the fact that the structure of zeros of certain Wronskians of orthogonal polynomials or
special functions has recently been studied numerically (cf.~\cite{Clarkson,Fedler,FR} and the references therein),
where it is shown that they have highly regular conf\/igurations in the complex plane.
It turns out that the zeros of Wronskians for certain multiple orthogonal polynomials produce intriguing pictures as
well, which might be of independent interest.

Throughout this section, we take $r=2$ and denote the multi-index $\vec{n}$ by $(n,m)\in\mathbb{N}^2$.
Unless otherwise stated, the path associated with the Wronskian~\eqref{def:wronskian} is chosen in such a~way that in
each step it is increased by one in the horizontal direction, i.e.,
\begin{gather*}
(n,m)\to (n+1,m) \to (n+2,m) \to \cdots.
\end{gather*}
Other choices of the paths show similar behavior of the zeros.
Clearly, the structure of the roots in the complex plane depends on $l$, $\vec{n}$, the path chosen and the values of the
parameters, therefore, it is dif\/ficult to be described completely.
We thus have chosen a~few illustrative examples from numerical experiments for multiple Hermite and Laguerre polynomials
to show the fascinating conf\/igurations.

\pagebreak

\subsection{Zeros of the Wronskians for multiple Hermite polynomials}

\begin{figure}[t]
 \centering
 \includegraphics[width=6.8cm]{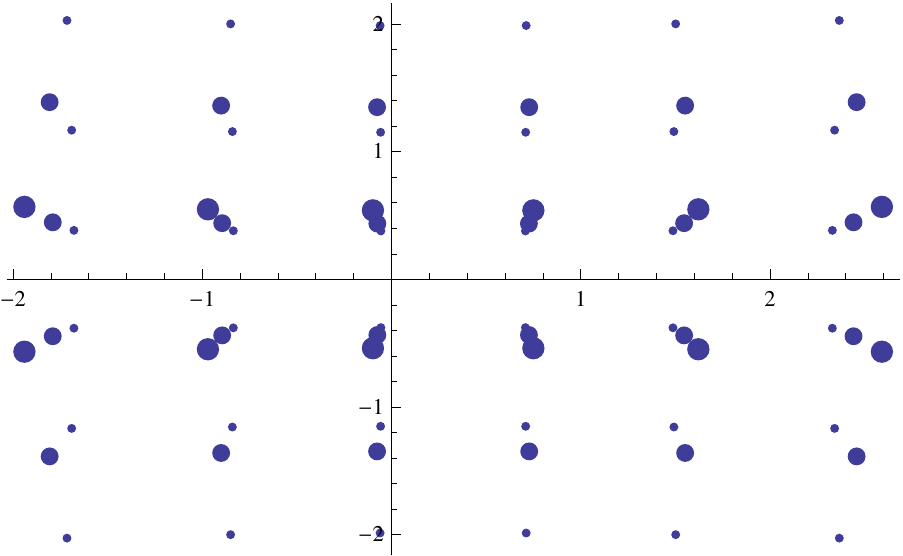} \qquad
 \includegraphics[width=6.8cm]{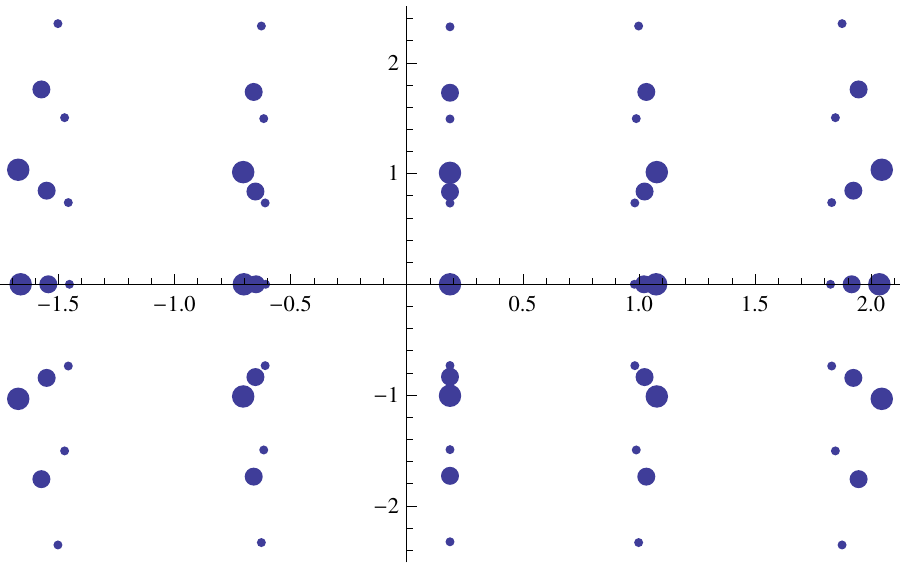}

  \caption{Zeros of the Wronskians for multiple Hermite polynomials with
$\vec{n}=(3,3)$, $\vec{c}=(1/3, 34/35)$ for $l=2,4,6$ (left) and $\vec{n}=(2,3)$, $\vec{c}=(1/3, 2/5)$ $l=3,5,7$
(right).
The size of points decreases as~$l$ increases.}\label{fig:1}
\end{figure}

\begin{figure}[t] \centering
\includegraphics[width=6.8cm]{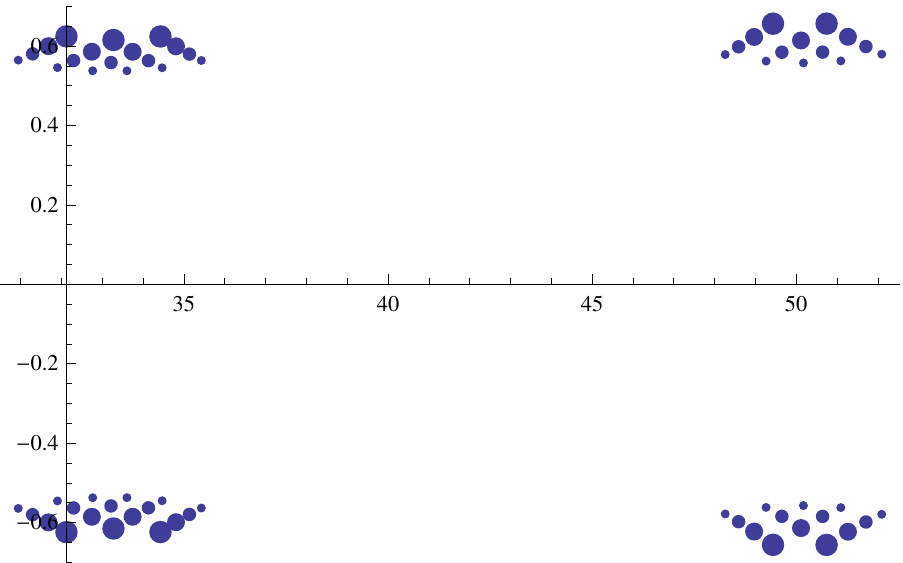}

 \caption{Zeros of the Wronskians for multiple Hermite
polynomials with $\vec{n}=(2,3),(3,4),(4,5),(5,6)$, $\vec{c}=(100,=200/3)$ for $l=2$.
The size of points decreases as $|\vec{n}|$ increases.}
\label{fig:2}
\end{figure}
\begin{figure}[t!]
\centering
 \includegraphics[width=6.8cm]{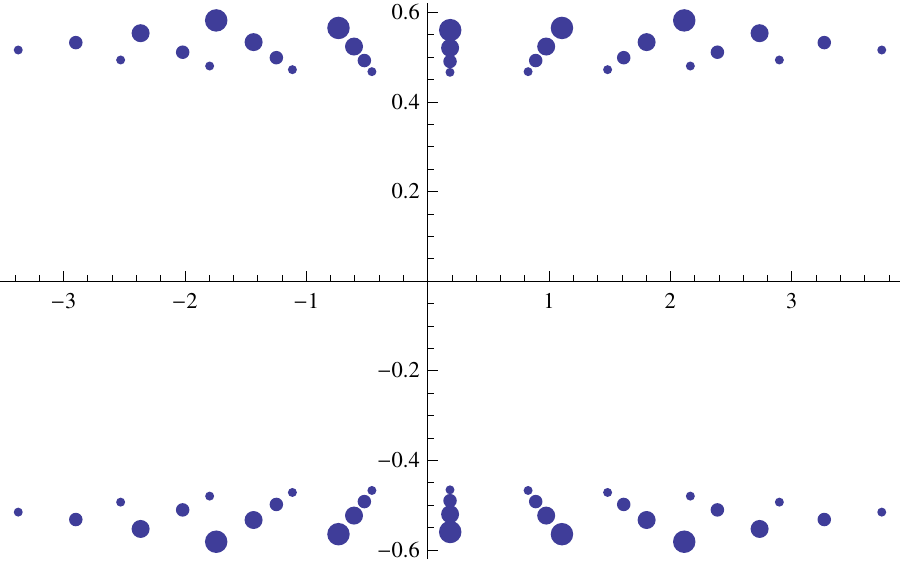} \qquad \includegraphics[width=6.8cm]{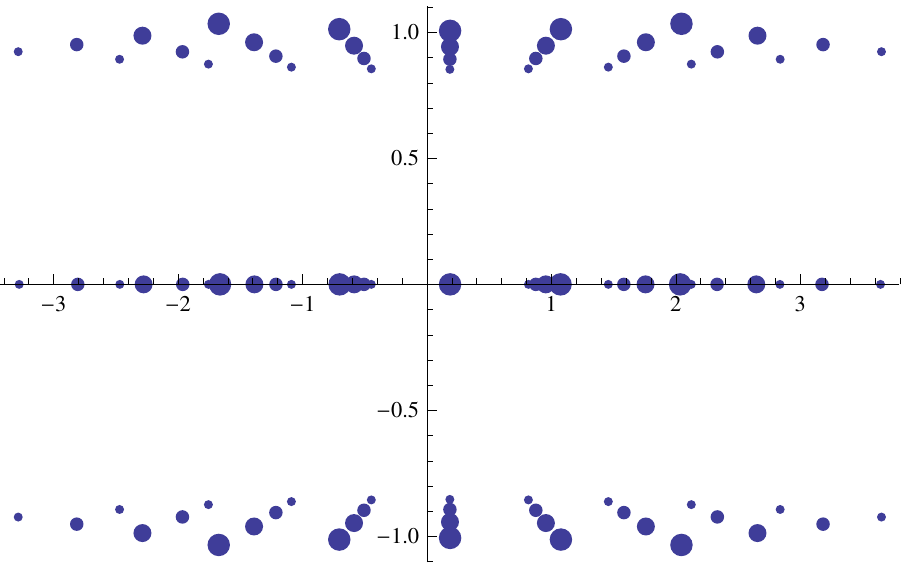}
\caption{Zeros of the Wronskians for multiple Hermite polynomials with $\vec{n}=(2,3),(3,4),(4,5),(5,6)$ and
$\vec{c}=(1/3,2/5)$ for $l=2$ (left) and $l=3$ (right).
The size of points decreases as $|\vec{n}|$ increases.}
\label{fig:3}
\end{figure}

The zeros of Wronskians for multiple Hermite polynomials numerically have roughly rectangular-like structure in the
complex plane.
In Fig.~\ref{fig:1} we plot zeros of the Wronskians for these polynomials by f\/ixing $\vec{n}$, the parameter $\vec{c}$
and increasing the length of the path (the size of points decreases as~$l$ increases).
We can see additional row of zeros on the real axis for~$l$ odd, as indicated by the f\/irst part of Theorem~\ref{thm:l
odd}.

If the two values in the parameter $\vec{c}$ dif\/fer too much, it seems that the zeros may separate into several
rectangles, which is shown in Fig.~\ref{fig:2} (the size of points decreases as $|\vec{n}|$ increase).
When~$l$ is odd, we have additional groups of zeros on the real line as projections of complex groups of roots.

The ef\/fect of increasing $|\vec{n}|$ is illustrated in Fig.~\ref{fig:3} for~$l$ even and odd respectively.
As $|\vec{n}|$ increases, the zeros are distributed in a~wider range.
Furthermore, if~$l$ increases, we can see more horizontal lines.

\begin{figure}[t]
\centering
\includegraphics[width=6.8cm]{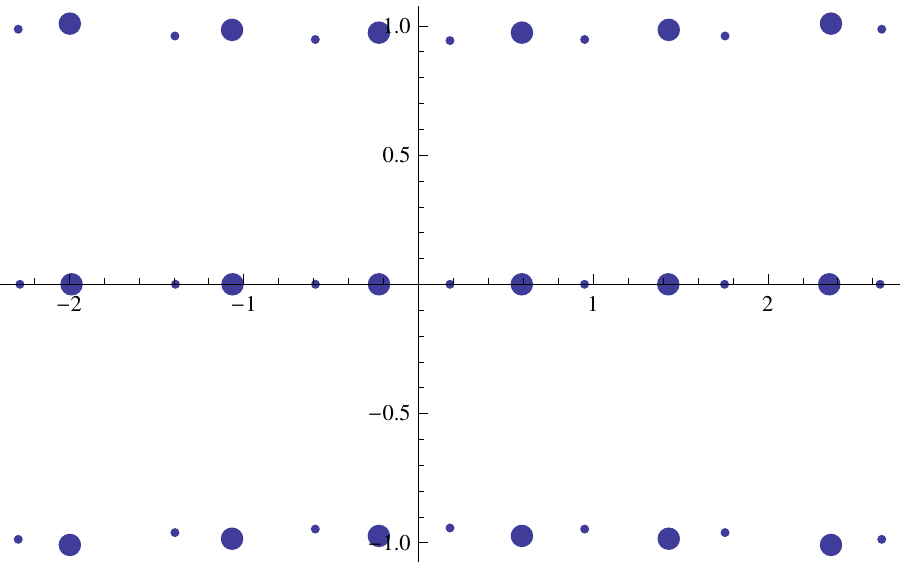} \qquad \includegraphics[width=6.8cm]{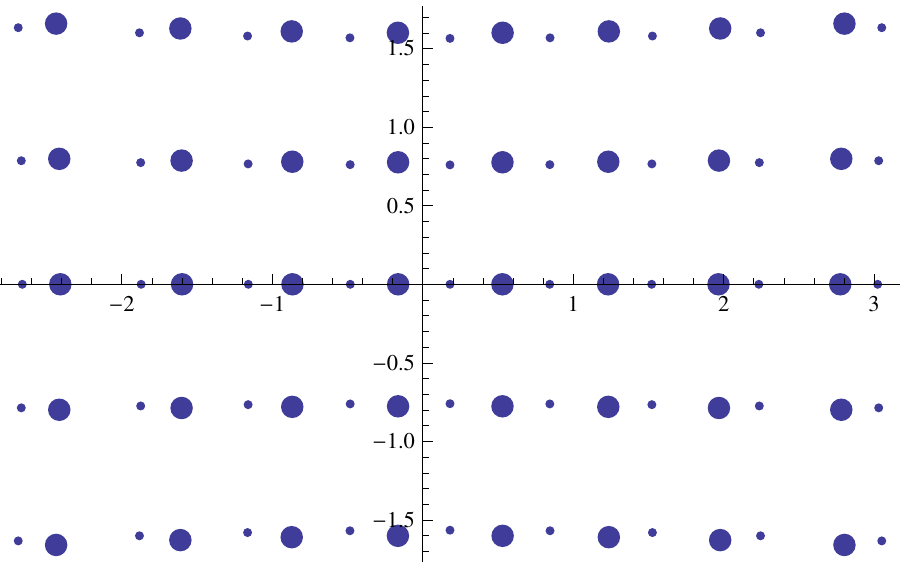}

\caption{Theorem~\ref{thm:l odd} for multiple Hermite polynomials with $\vec{c}=(1/3,2/5)$, $\vec{n}=(3,3)$ for $l=3$
(left) and $\vec{n}=(4,4)$ for $l=5$ (right).
The size of points decreases as $|\vec{n}|$ increases.}
\label{fig:4}
\end{figure}

We illustrate Theorem~\ref{thm:l odd} in Fig.~\ref{fig:4}.
We clearly see the interlacing of real zeros and regular conf\/igurations of zeros in the complex plane.
It seems that the interlacing property also appears on the other lines parallel to the real axis.
For even~$l$ the structure of complex roots is similar (there are no real roots in this case).

\begin{figure}[t]\centering
\includegraphics[width=6.8cm]{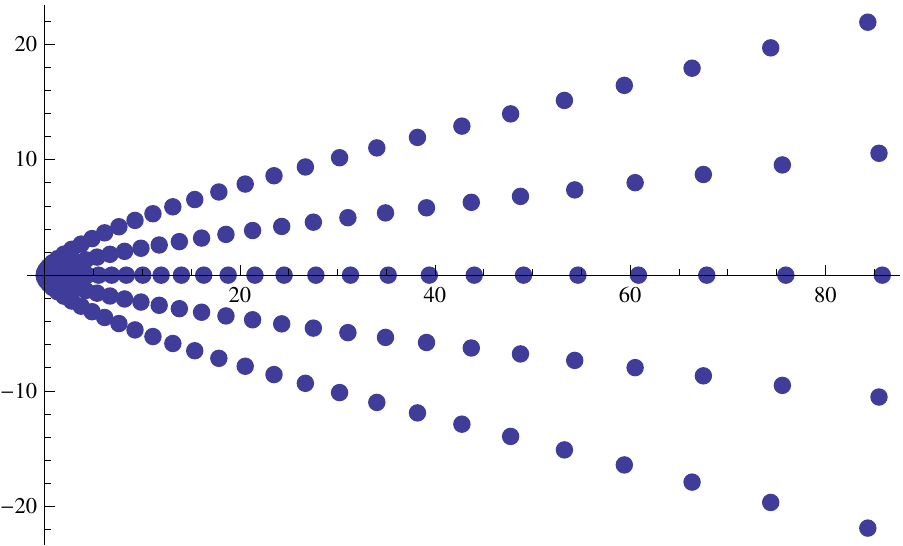} \qquad \includegraphics[width=6.8cm]{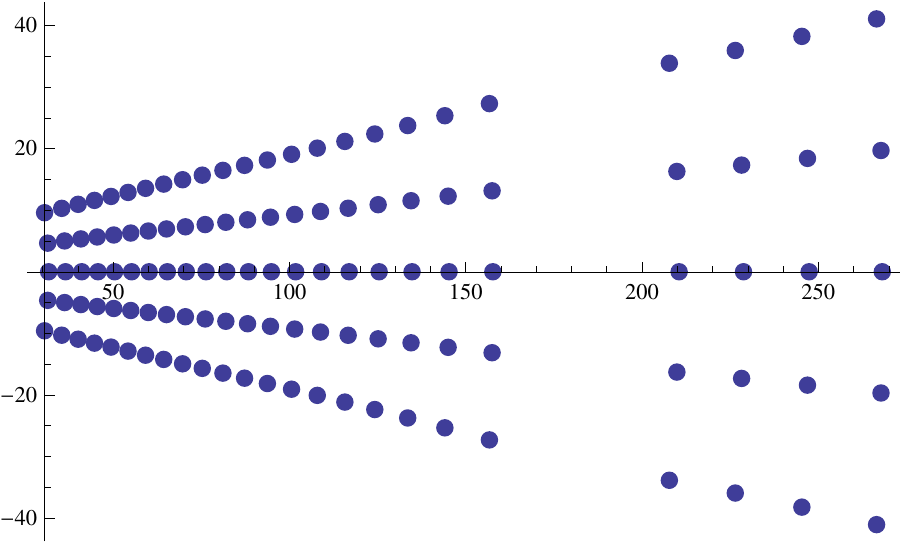}
\caption{Zeros of the Wronskians for multiple Laguerre polynomials of the f\/irst kind with $\vec{\alpha}=(1/2,1/3)$,
$\vec{n}=(10,20)$ (left) and $\vec{\alpha}=(200, 200/3)$, $\vec{n}=(4,20)$ (right) for $l=5$.}
\label{fig:5}
\end{figure}

\begin{figure}[t!]\centering
\includegraphics[width=6.8cm]{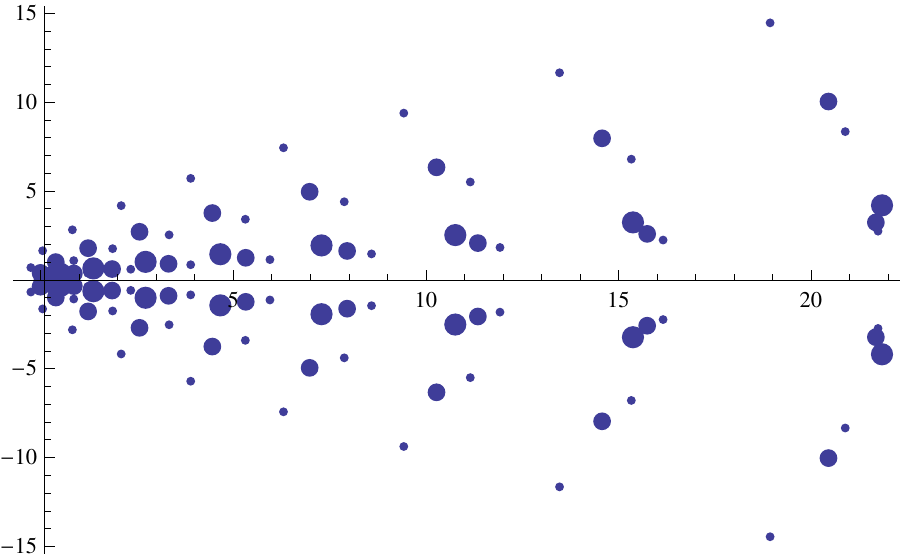}  \caption{Zeros of the Wronskians for multiple Laguerre polynomials of
the f\/irst kind with $\vec{n}=(4,5)$, $\vec{\alpha}=(1/2,1/3)$ for $l=2,4,6$.
The size of points decreases as~$l$ increases.}
\label{fig:6}
\end{figure}

\begin{figure}[t!]\centering
\includegraphics[width=6.8cm]{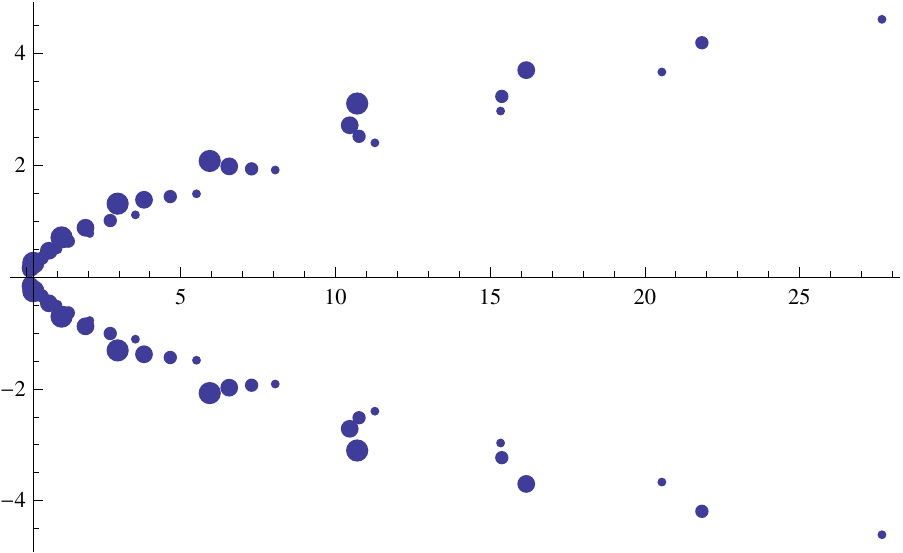} \qquad \includegraphics[width=6.8cm]{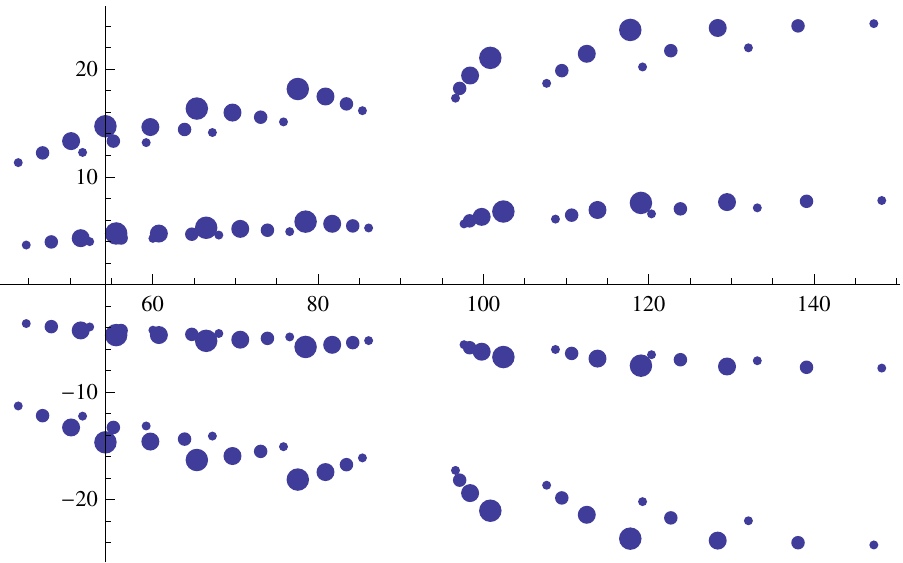}
\caption{Zeros of the Wronskians for multiple Laguerre polynomials of the f\/irst kind with $\vec{n}
=(2,3),(3,4),(4,5),(5,6)$, $\vec{\alpha}=(1/2,1/3)$ for $l=2$ (left) and $\vec{\alpha}=(100,200/3)$ for $l=4$ (right).
The size of points decreases as $|\vec{n}|$ increases.}
\label{fig:7}
\end{figure}

\begin{figure}[t!]\centering
\includegraphics[width=6.8cm]{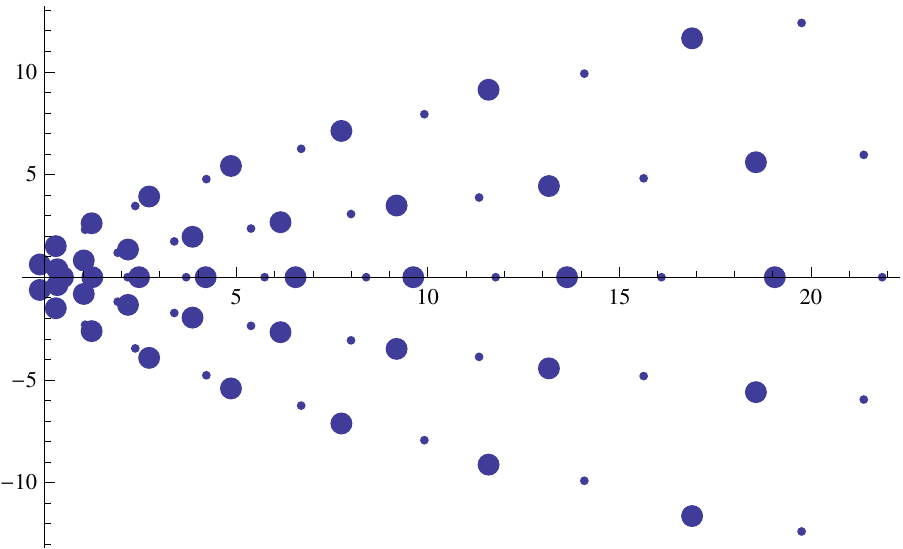} \qquad \includegraphics[width=6.8cm]{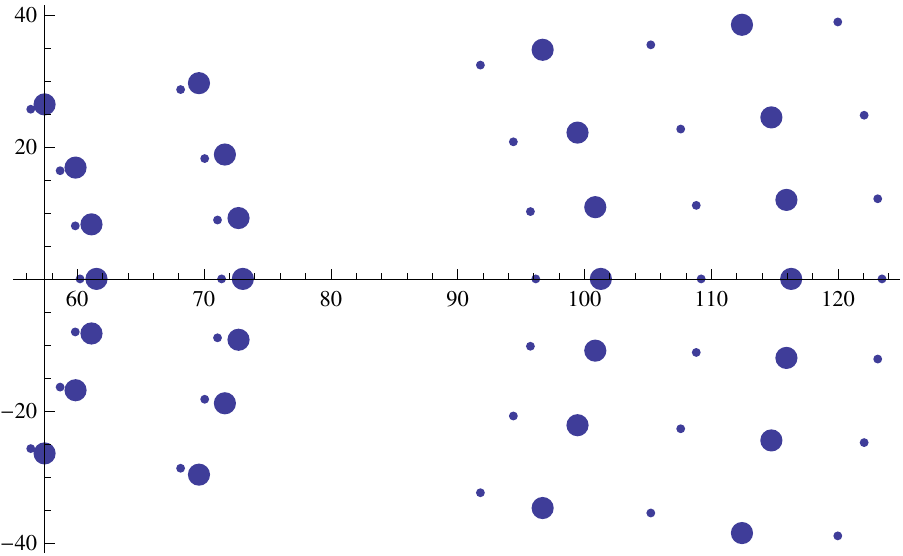}
\caption{Theorem~\ref{thm:l odd} for multiple Laguerre polynomials of the f\/irst kind with $\vec{n}=(4,4)$,
$\vec{\alpha}=(1/2,1/3)$ for $l=5$ (left) and $\vec{n}=(2,2)$, $\vec{\alpha}=(100,200/3)$ for $l=7$ (right).
The size of points decreases as $|\vec{n}|$ increases.}
\label{fig:8}
\end{figure}

\begin{figure}[t!]\centering
\includegraphics[width=6.8cm]{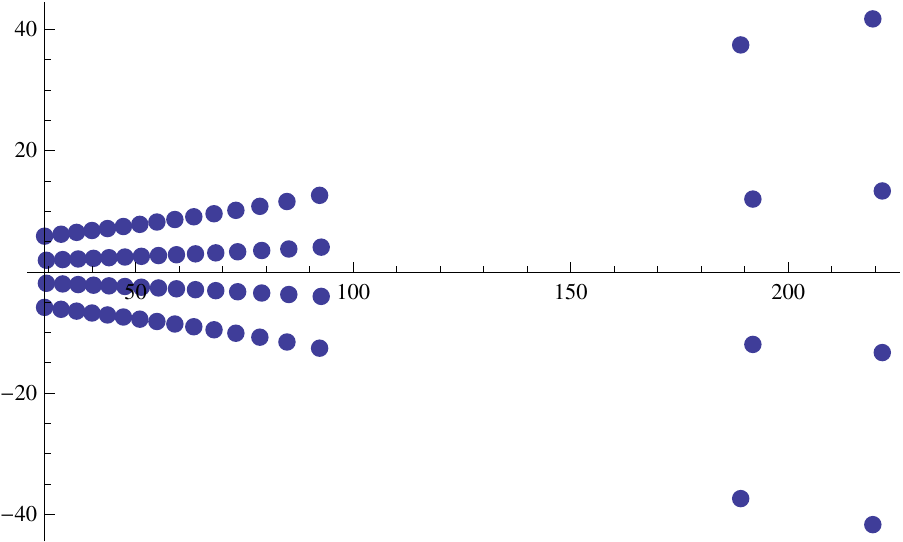} \caption{Zeros of the Wronskians for multiple Laguerre polynomials of
the second kind with $\vec{n}=(15,2)$, $\alpha=100$, $\vec{c}=(2,3/5)$ for $l=4$.}
\label{fig:9}
\end{figure}

\begin{figure}[t!]\centering
\includegraphics[width=6.8cm]{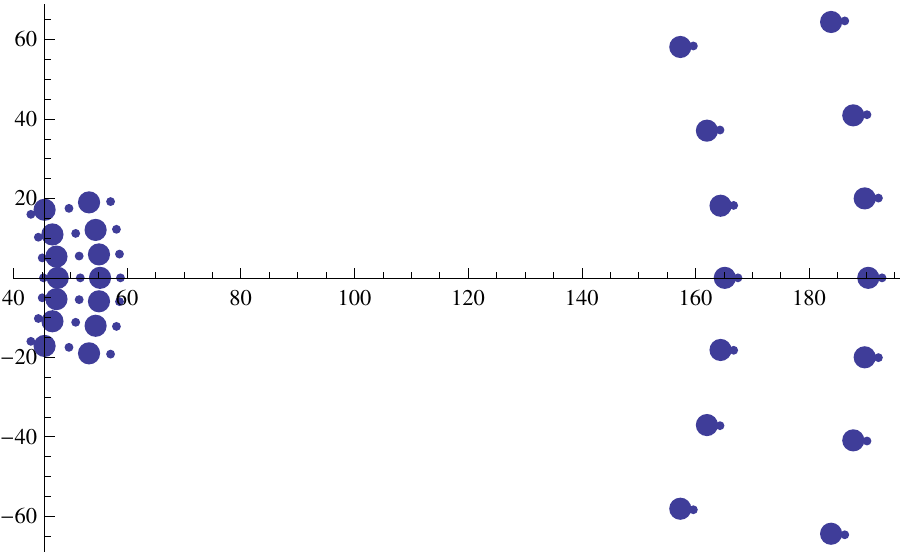} \caption{Theorem~\ref{thm:l odd} for multiple Laguerre polynomials of
the second kind with $\vec{n}=(2,2)$, $\alpha=100$, $\vec{c}=(2,3/5)$ for $l=7$.
The size of points decreases as $|\vec{n}|$ increases.}
\label{fig:10}
\end{figure}

\subsection{Zeros of the Wronskians for multiple Laguerre polynomials}

Since the observed structure of roots of the Wronskians for multiple Laguerre polynomials of the f\/irst and second kind
is numerically quite similar, we shall concentrate more on the multiple Laguerre polynomials of the f\/irst kind.
The conf\/igurations of zeros of Wronskians for multiple Laguerre polynomials of the f\/irst kind resemble (several)
parabolas (with additional zeros on the real line in case~$l$ is odd).
Sometimes it looks like zeros lie on arcs of circles with increasing radius.
As we change the parameters, we can observe that the zeros on the left can do not accumulate and, as in the case of
multiple Hermite polynomials, they can be grouped into several clusters with a~certain gap between them.
Fig.~\ref{fig:5} illustrates these observations.

The zeros of several Wronskians, if plotted together, also present nice interlacing properties.
Following the same strategy in previous section, we f\/ix $\vec{n}$, the parameter $\vec{\alpha}$ and increase~$l$ in
Fig.~\ref{fig:6} (compare with Fig.~\ref{fig:1}), while in Fig.~\ref{fig:7} we increase $\vec{n}$ and f\/ix other
parameters (compare with Fig.~\ref{fig:3}).
Theorem~\ref{thm:l odd} in this case is illustrated in Fig.~\ref{fig:8}.

Finally, the roots of the Wronskians for multiple Laguerre polynomials of the second kind are depicted in
Fig.~\ref{fig:9} and Theorem~\ref{thm:l odd} is illustrated in Fig.~\ref{fig:10}.

More pictures with dif\/ferent conf\/igurations of roots can be found in Mathematica f\/iles on the web-pages of the authors
(or available on request).
We believe the nice and regular geometric conf\/igurations generated from the zeros of Wronskians deserve further analytic
investigations.

\subsection*{Acknowledgements}

We thank the referees for helpful comments, suggestions, and pointing out the additional re\-fe\-rences~\cite{Elbert2,
Elbert1, Laforgia, Lorch}.
LZ is partially supported by The Program for Professor of Special Appointment (Eastern Scholar) at Shanghai Institutions
of Higher Learning (No.~SHH1411007) and by Grant SGST 12DZ 2272800 from Fudan University.
GF is supported by the MNiSzW Iuventus Plus grant Nr 0124/IP3/2011/71.

\pdfbookmark[1]{References}{ref}
\LastPageEnding

\end{document}